\documentclass[11pt]{article}

\input{epsf}

\usepackage{amsfonts}
\usepackage{amsmath}
\usepackage{amssymb}
\usepackage{pdfpages}

\newcommand{\Z}{\mathbb{Z}}
\newcommand{\R}{\mathbb{R}}

\newcommand{\rt}{\vartriangleright}

\newtheorem{theorem}{Theorem}[section]

\newtheorem{lemma}[theorem]{Lemma}
\newtheorem{proposition}[theorem]{Proposition}

\newtheorem{definition}[theorem]{Definition}

\newtheorem{example}[theorem]{Example}

\date{}

\begin{document}

\title{ Distributivity in Quandles and Quasigroups 
  }

\author{ 
Mohamed Elhamdadi \footnote{Email: \texttt{emohamed@math.usf.edu}}
\\ University of South Florida\\\\
{\it Dedicated to the memory of Jean-Louis Loday}
}

\maketitle

\vspace{10mm}

\begin{abstract}
Distributivity in algebraic structures appeared in many contexts  such as in quasigroup theory, semigroup theory and algebraic knot theory.  In this paper we give a survey of distributivity in quasigroup theory and in quandle theory.
\end{abstract}


\section{Introduction}

Quandles are in general non-associative structures whose axioms correspond to the algebraic distillation of the three Reidemeister moves in knot theory. Quandles appeared in the literature with many different names.  If one restricts himself to the most important axiom of a quandle which is the self-distributivity axiom (see definition below), then one can trace this back to 1880 in the work of  Pierce \cite{Peirce} where one can read the following comments { \it {"These are other cases of the distributive principle ....These formulae, which have hitherto escaped notice, are not without interest.", }}  Another early work fully devoted to self-distributivity appeared in 1929 by Burstin and Mayer \cite{BM} where normal subquasigroups are studied and an attempt is made to show
that every minimal subquasigroup of a finite distributive quasigroup is normal.  This is considered as the starting point for the investigation of normality problems in distributive
quasigroups.  In 1942 Mituhisa Takasaki \cite{Takasaki} introduced   the notion of kei (involutive quandle in Joyce's terminology \cite{Joyce}) as an abstraction of the notion of symmetric transformation.
The earliest known work on racks (see definition below) is contained
in the 1959 correspondence between John Conway and Gavin Wraith who studied
racks in the context of the conjugation operation in a group.   Around 1982,  Joyce \cite{Joyce} (used the term quandle) and Matveev \cite{Matveev} (who call them distributive groupoids)  introduced independently the notion of a quandle.   
Joyce and Matveev associated to each oriented knot $K$ a quandle $Q(K)$ called the knot quandle. The knot quandle is a complete invariant up to orientation.  Since then quandles and racks have been investigated by topologists  in order 
to construct knot and link invariants and their higher analogues (see for 
	example \cite{CKS} and references therein).  In 1986, Brieskorn \cite{Brieskorn} introduced the concept of automorphic set to describe a set $\Delta$ with a binary operation $*$ such that all  left multiplications $b \mapsto a*b$ are automorphisms of $\Delta$.  He considered the action of the braid group $B_n$   on the Cartesian product $\Delta ^n$ and introduced invariants of the orbit, for example, monodromy groups.  In 1991, Kauffman intoduced a similar notion called crystal (\cite{K&P} p. 186) as a generalization of the fundamental group of a knot in the sense that the crystal has more information than the fundamental group alone.  In 1992, Fenn and Rourke \cite{FR} showed that any codimension-two link has a fundamental rack which contains more information than the fundamental group.  They gave some examples of computable link invariants derived from the fundamental rack and explained the connection of the theory of racks with that of braids.  In 2003, Fenn, Rourke and Sanderson  \cite{FRS}  introduced rack homology.  This (co)homology was modified in 1999 by Carter et al. \cite{CJKLS} to give a cohomology theory for quandles.   This cohomology was used to define state-sum invariant for knots in three space and knotted surfaces in four space.  A nice survey paper on quandle ideas is a paper by Scott Carter \cite{Carter} showing the applications of quandle cocycle invariants.

\noindent
In this paper, we give a survey of distributivity in quasigroup theory and in quandle theory. \\
  \noindent 
 In Section 2, we review the basics of quandles and give examples.  Section 3 deals with the problem of classification of quandles.  In section 4 we relate quandles to quasigroups and Moufang loops.  Section 5 deals with  the quandle cohomology and cocycle knot invariants. 
 \noindent

\section{ Basics of quandles} 
We start by reviewing the basics of quandles and give some examples.

\begin{figure}\label{Rmoves}
\begin{center}
\mbox{
\includegraphics[scale=0.55]{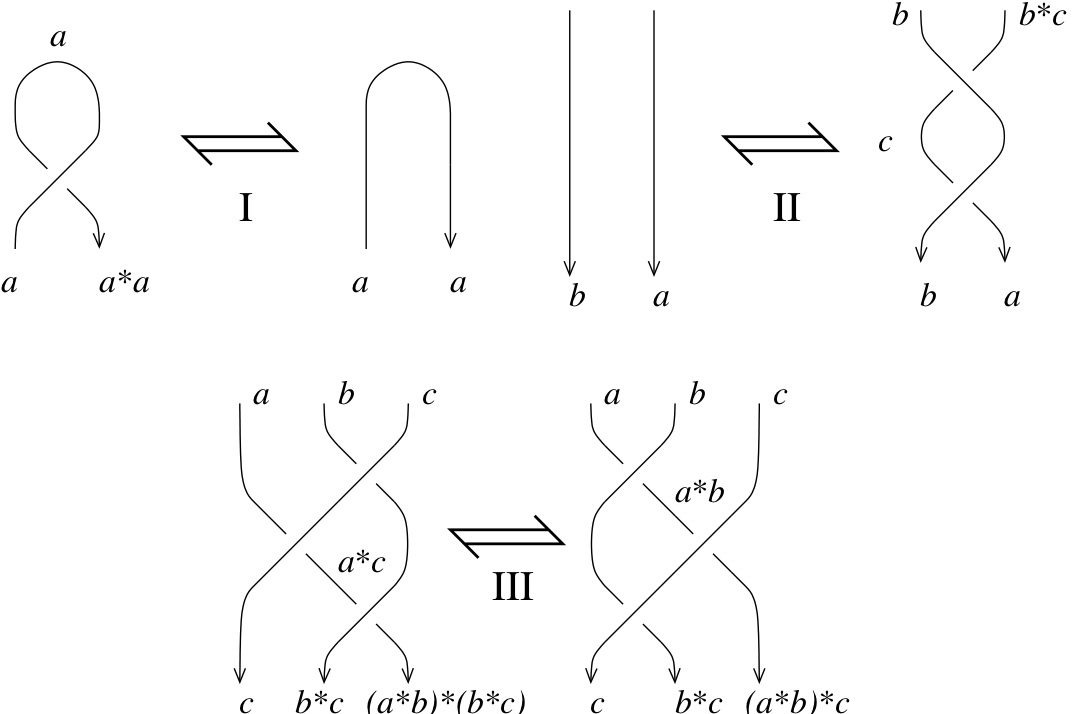}
}
\end{center}
\caption{ Reidemeister moves and the quandle identities  }
\end{figure}

\begin{definition} {\rm \cite{Joyce}\label{Joyce}
A {\it quandle}, $X$, is a set with a binary operation $(a, b) \mapsto a * b$
such that

(1) For any $a \in X$,
$a* a =a$.

(2) For any $a,b \in X$, there is a unique $c \in X$ such that 
$a= c*b$.

(3) 
For any $a,b,c \in X$, we have
$ (a*b)*c=(a*c)*(b*c). $
}
\end{definition}
 \noindent
Axiom (2) states that for each $u \in X$, the map $R_u:X \rightarrow X$ with $R_u(x):=x*u$ is a bijection.  
\noindent
 The axioms for a quandle correspond respectively to the 
Reidemeister moves of type I, II, and III as can be seen from Figure \ref{Rmoves}.  Quandles have been used to study colorings of knots and links and to define some of their invariants, see for example \cite{CES}.   

\bigskip

\noindent
Here are some examples of quandles: \\
$\bullet$ Any set $X$ with the operation $x*y=x$ for all $x,y \in X,$ is
a quandle called the {\it trivial} quandle.
\\ $\bullet$ Any group $ X = G$ with conjugation $ a * b = bab^{-1}$ is a quandle.
\\ $\bullet$  Let $n$ be a positive integer.
For elements  
$i, j \in \mathbb{Z}_n$ (integers modulo $n$), 
define
$i\ast j \equiv 2j-i \pmod{n}$.
Then $\ast$ defines a quandle
structure  called the {\it dihedral quandle},
  $R_n$.
This set can be identified with  the
set of reflections of a regular $n$-gon
  with conjugation
as the quandle operation.  If we denote the group of symmetry of a regular  $n$-gon by $ D_n = < u,v \; | \;u^n=1, v^2=1, vuv=u^{-1}>$, then conjugation on reflections is given by $ (u^iv)*(u^jv) = u^jvu^iv(u^jv)^{-1}  = u^ju^{-i}vu^{-j} = u^{2j-i}v.$
\\ $\bullet$ 
A group $X=G$ with operation $x*y=yx^{-1}y$ is called the $core$ quandle of $G$, denoted $Core(G)$.  
\\ $\bullet$  For any abelian group $M$ and automorphism $t$ of $M$ define a quandle structure on $M$ by $x*y = t(x-y)+y$.  This is called an {\it Alexander quandle}.
\\ $\bullet$ A generalization of the last example is, let $G$ be a group and $\phi$ be an automorphism of $G$, then define a quandle structure on $G$ by $x*y=\phi(xy^{-1})y$.  Further, Let $H$ be a subgroup of $G$ such that $\phi(h)=h,$ for all $ h \in H$.  Then $G/H$ is a a quandle with operation $Hx *Hy=  H \phi(xy^{-1})y$.  It is called the $homogeneous$ quandle $(G,H,\phi)$.
\\ $\bullet$ Let $<\;,\;>:\R^n \times \R^n \rightarrow \R$ be a symmetric bilinear form on $\R^n$.  Let $X$ be the subset of $\R^n$ consisting of vectors $x$ such that $<x,x>\neq 0$.  Then the operation $$x*y=\frac{2<x,y>}{<x,x>}y-x$$ defines a quandle structure on $X$.  Note that, $x*y$ is the image of $x$ under the reflection in $y$.  This quandle is called a $Coxeter$ quandle.  

\noindent
A function $\phi: (X, *) \rightarrow  (Y,\rt)$ is a quandle {\it homomorphism}
if $\phi(a \ast b) = \phi(a) \rt \phi(b)$ 
for any $a, b \in X$.   Axiom (3) of the definition \ref{Joyce} state that for each $u \in X$, the map $R_u$  
is a quandle homomorphism.   Let {\rm Aut(X)} denotes the automorphism group of $X$.   The subgroup of {\rm Aut(X)}, generated by the permutations $R_x$, is called the $inner$ automorphism group of $X$ and denoted by {\rm Inn}$(X)$.   By axiom (3) of the definition \ref{Joyce}, the map $R:X \rightarrow {\rm Inn}(X)$,  sending $u$ to $R_u$, satisfies  the equation $R_z\;R_y=R_{y*z}\;R_z, \;\;$ for all $ y,z \in X$, which can be written as $R_z\;R_y\;{R_z}^{-1}=R_{y*z}.$ Thus, if the group  {\rm Inn}$(X)$ is considered as a quandle with conjugation then the map $R$ becomes a quandle homomorphism.  The subgroup of {\rm Aut(X)}, generated by  $R_xR_y^{-1},$ for all $ x, y \in X$, is called the $transvection$ group of $X$ denoted by $Transv(X)$.  It is well known (see for example \cite{Joyce}) that  the $transvection$ group is a normal subgroup of the inner group and the later group is normal subgroup of the   automorphism group of $X$.  The quotient group ${\rm Inn(X)}/ {\rm Transv}(X)$ is a cyclic group (see  \cite{Joyce}).   For each $u \in X$, let denotes the left multiplication by the map $L_u:X \rightarrow X$ with $L_u(x):=u*x$.  We list some properties and some definitions  of quandles below.

\begin{itemize}
\setlength{\itemsep}{-3pt}
\item
A quandle $X$ is {\it involutory}, or a {\it kei}, 
if the right translations  are  involutions: ${ R}_a^2 ={\rm  id},$  for all $a \in X$.

\item

A quandle is {\it faithful} if the mapping $a \mapsto {R}_a$ is 
an injection 
from $X$ to ${\rm Inn}(X)$.

\item
A quandle is {\it connected} if ${\rm Inn}(X)$ acts transitively on $X$.

\item
A {\it Latin quandle} is a quandle such that for each $a \in X$, the left translation ${ L}_a$ is a bijection. 
That is, the multiplication table of the quandle is a Latin square. 

\item
A quandle $X$ is {\it medial} if 
$(a*b)*(c*d)=(a*c)*(b*d)$  for all  $a,b,c,d \in X$. 
It is well known that a quandle is medial iff its tranvection group is abelian, that is why it is also called {\it abelian}. 
It is known and easily seen that 
every Alexander quandle is medial.

\item

A quandle $X$ is called $simple$ if the only surjective quandle homomorphisms
on $X$ have trivial image or are bijective.

\end{itemize}

\section{The problem of classification of quandles}
The problem of classification of quandles and racks was attempted by many authors mainly because computable invariants  of knots such as, for example,  the quandle cocycle invariant of Carter et al. \cite{CJKLS, CKS},  and enhancement of counting homomorphisms from the knot quandle to a fixed quandle of Nelson et al. \cite{NN, N} can be defined from quandles.  Racks and quandles are used in the classification of pointed Hopf algebras \cite{AG} since they help in the understanding of Yetter-Drinfeld modules over groups.  Below, we give a survey of the classification of finite quandles.

  In 2003, Nelson gave a classification of finite Alexander quandles  proving the following
  \begin{theorem}   \cite{N1}
  Two finite Alexander quandles $M$ and $N$ of the same cardinality are isomorphic as quandles if and only if $(1-t)M$ and $(1-t)N$ are isomorphic as $\Z[t,t^{-1}]$-modules. 
  \end{theorem}
    As a consequence of this theorem Ho and Nelson \cite{HN}  computed  isomorphism  classes of quandles up to order $5$ and their automorphim groups.   Quandles of order 6, 7 and 8 were given by Henderson, Macedo and Nelson in \cite{ HMN} but isomorphism class representatives were not determined.   In 2006, Nelson and Wong \cite{NW} obtained the orbit decomposition of finite quandles:  A subset $A$ of a quandle $X$ is said to be $X$-complemented if the complement of $A$ in $X$ is a subquandle of $X$.  They proved the following 
    \begin{theorem}  \cite{NW}  
     Up to isomorphism, every finite quandle has a unique decomposition into subquandles $A_1, A_2,  \ldots ,  A_n$ such that every $A_j$ is $X$-complemented and no proper subquandle of any $A_j$ is $X$-complemented.
      
  \end{theorem}
Independently around the same time Yetter et al.  \cite{EGTY} obtained a similar decomposition theorem for quandles in terms of an operation of "semidisjoint union'', showing that all finite quandles canonically decompose via iterated semidisjoint unions into connected subquandles.  Murillo and Nelson  \cite{MN} proved in 2006 that there are 24 isomorphism classes of Alexander quandles of order 16. 
In  \cite{EMR} quandles up to order 9 were classified, automorphism groups of quandles (with orders up to $7$) were determined and the automorphism group of the dihedral quandle $R_n$ was proven to be isomorphic to the affine group of $\Z_n$ .  The
number of isomorphism of quandles of order 3, 4, 5, 6, 7, 8 and 9  are
respectively 3, 7, 22, 73, 298, 1581, 11079.   The list of isomorphism classes can be found in https://sites.google.com/a/exactas.udea.edu.co/restrepo/quandles.  Independently the same classification result was obtained in \cite{OEIS} by McCarron.\\ 
In \cite{Hou}  it was shown first that the isomorphism class of an Alexander quandle $(M,*)$ is determined by the isomorphism type of the $\Lambda$-module $(1-t)M$ and the cardinality of the quotient $A/K,$ where $A$ is the annihilator of $(1-t)$ in $M$, $K=A\cap(1-t)M$ and $\Lambda=\Z [t,t^{-1}]$. This recovers a result of Sam Nelson \cite{N1}.  The structure of the automorphism group of a general Alexander quandle $(M,*)$ is completely determined (see \cite{Hou} for more details).  Enumeration of Alexander quandles has been much improved.  Edwin Clark computed the number of Alexander quandles of orders up to $255$ (see http://oeis.org/A193024, for more details) based on results from \cite{Hou1} which
contains other interesting enumeration results concerning Alexander
quandles.  More sequences related to quandles can be found on http://oeis.org.  

\begin{example}{\rm
One way of describing a finite quandle is by the Cayley table.  Since by the second axiom of a quandle right multiplication by a fixed $i$, $R_i: j \mapsto j*i$  is  is a permutation.  We then can describe each quandle by writing each column $R_i$ of the Cayley table as a product of disjoint cycles.  Here we include the list of quandles of order 4. The notation $(1)$ in the table means that the permutation is the identity permutation.  For example the quandle $Q_5$ is the set $\{1,2,3,4\} $ where $R_1$ is the identity permutation, $R_2$ is the transposition sending $3$ to $4$, $R_3$ is the transposition sending $2$ to $4$ and $R_4$ is the transposition sending $2$ to $3$.
\small{
\begin{table}
\begin{center}
\begin{tabular}{|l|c|c||} \hline 
Quandle & Disjoint Cycle Notation for the Columns of the Quandle  
\\ \hline
$Q_1$ & $(1),(1),(1),(1)$   \\ \hline
$Q_2$ & $(1),(1),(1),(23)$   \\ \hline
$Q_3$ & $(1),(1),(1),(123)$   \\ \hline
$Q_4$ & $(1),(1),(12),(12)$   \\ \hline
$Q_5$ & $(1),(34),(24),(23)$   \\ \hline
$Q_6$ & $(34),(34),(12),(12)$   \\ \hline
$Q_7$ & $(234),(143),(124),(132)$   \\ \hline
\end{tabular}\end{center}
\caption{Quandles of Order 4 in Terms of Disjoint Cycles of Columns}
\label{qmoduletable3}
\end{table}
}
}   
\end{example}
  Using computers the search space in general becomes too large to obtain the computation of all quandles up to isomorphism for higher cardinality. Clearly, this depends on the algorithm used to find quandles.   However if one restricts himself to the subclass of connected quandles then classification becomes more accessible
to calculation in somehow comparable way to the classification of finite groups.  In \cite{Clauw}, Clauwens studied connected quandles and proved the following
\begin{proposition}{\rm \cite{Clauw}
If $f:Q \rightarrow P$ is a surjective quandle homomorphism and $P$ is connected then for all $x, y \in P$, there is a bijection between $f^{-1}(x)$ and $f^{-1}(y)$.  In particular  the cardinality of $P$ divides the cardinality of $Q$.
}
\end{proposition}
This allowed him to obtain isomorphism classes of connected quandles up to order 14, in
particular he showed that there is no connected quandle of order 14. In \cite{Ven}, Vendramin extended Clauwens results to the list of all connected quandles of orders less than 36.  He used the classification of transitive groups and the program
described in \cite{EGTY} based mainly on the following 

\begin{theorem} \cite{Ven}
Let $X$ be a connected quandle of cardinality $n$.  Let $x_0\in X$ and $z=R_{x_0}$ be the right multiplication by $x_0$, $G={\rm Inn}(X)$ and $H=Stab_G(x_0)=\{g \in G, gx_0=x_0\}$.  Then (1) $G$ is a transitive group of order $n$, (2) $z$ is central element of $H$ and (3) $X$ is isomorphic to the homogeneous quandle $(G,H,I_z),$ where  $I_z$ is the conjugation by $z$.

\end{theorem}
A complete list of isomorphism classes of quandles with up to 6 elements appeared in the appendix \cite{CKS} .

\section{Quandles and quasigroups} 
 In this section we will discuss the relation between left and right distributive quasigroups and the following types of quandles: Alexander, Latin and medial quandles.  Two connections between quasigroups and quandles were established in \cite{Smith1}. \\
Self-distributivity appeared in 1929 by Burstin and Mayer \cite{BM} where they studied quasigroups which are left- and
right-distributive.  They stated that there are none of orders 2 and 6,
observed that the group of automorphisms is transitive, and showed that
such a quasigroup is idempotent. 

 \begin{definition}{ \rm \cite{Bruck}
(1) A quasigroup is a set $Q$ with a binary operation $*$ such for all $u \in Q$ the right translation $R_u$  and left translation $L_u$ by $u$ are both permutations. \\
(2) If the operation $*$ has an identity element $e$ in $Q$  then the quasigroup is called a $loop$ and denoted $(Q,*,e)$. 
}
 \end{definition}  
\noindent
 Quasigroups differ from groups in the sense that they satisfy identities which usually conflict with associativity.  Distributive quasigroups have transitive groups of automorphisms but the only group with this
property is the trivial group.  In \cite{Stein} it is shown that there are no right-distributive quasigroups whose order is twice an odd number.   Right-distributive quasigroups are intimately connected with the binary operation of a conjugation in a group since in a right-distributive quasigroup it holds that $R_{y*z}=R_zR_yR_z^{-1}$ and the mapping $x \mapsto R_x$ is injective. We will see below that distributive quasigroups relate to Moufang loops.
\begin{definition} { \rm \cite{Bruck}
Let $(M,*)$ be a set  with a binary operation.  It  is called a \emph{Moufang loop} if it is a loop such that the binary operation satisfies one of the following equivalent identities:
\begin{eqnarray}
 x*(y*(x*z))=((x*y)*x)*z, \label{mouf1}\\
 z*(x*(y*x))=((z*x)*y)*x, \label{mouf2}\\
 (x*y)*(z*x)=(x*(y*z))*x. \label{mouf3}
\end{eqnarray}
}
\end{definition}
As the name suggests, the Moufang identity is named for Ruth Moufang who discovered
it in some geometrical investigations in the first half of this century  \cite{Moufang}.  Moufang loops differ from groups in that they need not be associative. A Moufang loop that is associative is a group. The Moufang identities may be viewed as weaker forms of associativity.  The typical examples include groups and the set of nonzero octonions which gives a nonassociative Moufang
loop.
\begin{theorem} (Moufang's Theorem)
Let $a, b, c$ be three elements in a commutative Moufang loop (abbreviated CML) $M$ for which the relation $(a*b)*c=a*(b*c)$ holds.  Then the subloop generated by them is associative and hence is an Abelian group.

\end{theorem}
A consequence of this theorem is that every two elements in CML generate an Abelian subgroup.  
\noindent
Let $(X,*)$ be a  right-distributive quasigroup.  Then $(x*x)*x=(x*x)*(x*x)$ which implies that each element is idempotent and $(X,*)$ is then a Latin quandle.  Fix $a \in X$ and define the following operation, denoted $+$, on $X$ by $x+y:={R_a}^{-1}(x)*{L_a}^{-1}(y)$.  Then $a+y=y$ and  $y+a=y$.  Thus $(X, +, a)$ is a loop.  Therefore any   right-distributive quasigroup satisfying one of the Moufang identities (\ref{mouf1}), (\ref{mouf2}) and (\ref{mouf3}) is a Moufang loop.  Remark that $R_a(x)+L_a(y)=x*y$.  The Moufang loop is commutative if and only if 
\begin{eqnarray}\label{medial}
(u*v)*(w*z)=(u*w)*(v*z)
\end{eqnarray}
A magma $(X,*)$ that satisfies equation (\ref{medial}) is said to be $medial$ (Belousov \cite{Bel2}) or $abelian$ (Joyce \cite{Joyce}).  
The Bruck-Toyoda theorem gives the following characterization of medial quasigroup.  Given an Abelian group $M$, two commuting automorphisms $f$ and $g$ of $M$ and a fixed element $a$ of $M$, define an operation $*$ on $M$ by $x*y=f(x)+g(y)+a$.  This quasigroup is called {\it affine} quasigroup.  It's clear that $(M,*)$ is a medial qasigroup.  The Bruck-Toyoda theorem states that every medial quasigroup is of this form, i.e. is isomorphic to a quasigroup defined from an abelian group in this way.  Belousov gave the connection between distributive quasigroups and Moufang loops in the following

\begin{theorem}   \cite{Bel2}
If $(X,*)$ be a distributive quasigroup then for all $a \in X$, $(X,+,a)$ is a commutative Moufang loop.

\end{theorem}
Now let  $(X,*)$ be a Latin quandle (that is right-distributive quasigroup), then the automorphism $\phi=R_a$ satisfies $2\phi(a)=a$.  If the order of $a$ is odd then one can write $\phi(a)=\frac{1}{2}a$.  The map $x\mapsto 2x$ being a homomorphism is equivalent to $(x+y)+(x+y)=(x+x)+(y+y),$ (mediality property).\\
Recall that a $magma$ is a set with a binary operation.  we have the following question: do the following three properties imply associativity for a finite magma $(X,+)$?

1. $(X,+)$ is a commutative loop with identity element  0.

2. For all $x,y$ in $X$ we have the identity
   $ (x+y) + (z + z) = (x + z) + (y + z).$

3. There is an automorphism $f$ of $(X,+)$ satisfying
    $f(x) + f(x) = x$ for all $x$.
 (in other words,  the map $x \mapsto 2x$ is onto and $(x + x) + (y + y) = (x+y) + (x + y).$

In fact, if $(X,+)$ is a loop satisfying condition 2, then $(X,+)$ is a commutative Moufang loop, necessarily satisfying the other conditions. There exist nonassociative commutative Moufang loops. The smallest order at which such loops occur is 81, and there are, in fact, two such loops of that order.  The easier to describe of the two commutative Moufang loops of order $81$ is the one of exponent $3$. Special thanks to Michael Kinyon and David Stanovsky for telling us about the following example and some other results about quasigroups.   Let $F= {\mathbb Z}_3$ and on $F^4$, define
$$(x_0,x_1,x_2,x_3) + (y_0,y_1,y_2,y_3) =
(x_0 + y_0 + (x_1 - y_1)(x_2 y_3 - x_3 y_2),x_1 + y_1,x_2 + y_2,x_3 + y_3),$$
This is very first known example, published by Bol, who attributed it to Zassenhaus \cite{Bol}.\\ 
The construction from loops to quandles 
requires the maps
$x\mapsto 2x$ to be a bijections as well as a homomorphisms. Is this guaranteed for commutative Moufang loops?  Every abelian group is a commutative Moufang loop, so squaring is not always a bijection, of course.  For the two examples we mentioned above (loops of order $81$),  the answer is yes.  Any commutative Moufang loop modulo its center will have exponent $3$. If you have a commutative Moufang loop which is indecomposable in the sense that it is not a direct product of smaller loops, then it will have order a power of 3. Nonassociativity starts showing up at order 81. Classification of commutative Moufang loops of higher order has not been worked out in detail because of the computational difficulties. Much literature has been about free commutative Moufang loops of exponent $3$, because they turn out to be finite and of order $3^n$.  Quandles which are also quasigroups correspond to a class of loops known as Bruck loops.
Commutative Moufang loops have been investigated in detail by Bruck and Salby
\begin{theorem}  \cite{Bruck}
If $(X,+)$  is a commutative Moufang loop then $X=A\times B$ is a direct product of an abelian group $A$ with order prime to $3$ and a commutative Moufang loop of order $3^k$.

\end{theorem}
Latin quandles are right distributive quasigroups and left-distributive Latin quandles are distributive quasigroups.  Belousov's theorem tells us that if $(X,*)$ is left-distributive Latin quandle then $(X,+)$ is a commutative Moufang loop and then Bruck-Slaby theorem tells us that $(X,*)$ is affine over a commutative Moufang loop, and then medial.
The smallest Latin quandle that is not left distributive is of order $15$ and was found by David Stanovsky (see \cite{Stan}, p 29) using an automatic model builder SEM for
all quasigroups satisfying left distributivity, but not mediality.  This motivated Jan Vlachy \cite{Vla} to look for a more theoretical argument that would explain
the nonexistence of any smaller quasigroups of this kind and proved that there are exactly two non-isomorphic types of these smallest non-right-distributive left-distributive quasigroups with
15 elements.   He constructed them
explicitly using the Galkin's representation \cite{galkin1}.  In the survey paper  \cite{galkin}, page $950$, Galkin states that nonmedial quasigroups of order less than 27 appear only in orders 15 and 21 and are given by the following construction: Define a binary operation on $\Z_3 \times \Z_p$ by 
$$(x, a)*(y, b)=(2y-x, -a + \mu(x-y) b + \tau(x-y) ) \quad x, y \in \Z_3, \ a, b \in \Z_p,$$ where $\mu(0)=2$, $\mu(1)=\mu(2)=-1$, and $\tau: \Z_3 \rightarrow \Z_p$ is such that $\tau(0)=0$. 
This construction was generalized by replacing  $\Z_p$ by any abelian group $A$ in \cite{CEHSY}.  Let  $A$ be an abelian group, also regarded naturally as a $\Z$-module. 
Let $\mu: \Z_3  \rightarrow \Z$ ,  $\tau: \Z_3 \rightarrow A$ be functions.   
These functions $\mu$ and $\tau$ need not 
be homomorphisms.
Define a binary operation on $\Z_3 \times A$ by 
$$(x, a)*(y, b)=(2y-x, -a + \mu(x-y) b + \tau(x-y) ) \quad x, y \in \Z_3, \ a, b \in A. $$
 
\begin{proposition} \label{galprop}\cite{CEHSY}
 For any abelian group $A$, 
the above operation $*$ defines a quandle structure on $\Z_3 \times A$ if 
$\mu(0)=2$, $\mu(1)=\mu(2)=-1$, and $\tau(0)=0$. 
\end{proposition}
This quandle $(\Z_3 \times A, *)$ is called the {\it Galkin quandle  }
and denoted by $G(A, \tau)$. 

\begin{lemma}\label{c1c2lem}  \cite{CEHSY}
For any abelian group $A$  and $c_1, c_2 \in A$, 
$G(A,  c_1, c_2)$ and $ G(A,  0, c_2 - c_1 ) $
are isomorphic. 
\end{lemma}
Various properties of Galkin quandles were studied in \cite{CEHSY} and their classification in terms of pointed
abelian groups was given.  We mention a few properties.  Each $G(A, c)$ is connected but not Latin unless $A$ has odd order, $G(A, c)$ is non-medial unless $3A = 0$   \\
We conclude with the folowing properties relating distributivity and mediality to quandles \cite{CEHSY}:  
Alexander quandles are left-distributive and medial. 
It is easy to check  that for a finite Alexander quandle $(M, T)$ with $T \in {\rm Aut}(M)$,
the following are equivalent:
(1) $(M,T)$ is connected, 
(2) $(1-T)$ is an automorphism of $M$, and 
(3)  $(M,T)$ is Latin. 
It was also proved by Toyoda~\cite{Toyo} that 
a Latin quandle is Alexander if and only if it is medial. 
As  noted by Galkin, $G(\Z_5,0)$ and $G(\Z_5,1)$ are  the smallest non-medial Latin quandles
and hence the smallest non-Alexander Latin quandles. \\
We note that medial 
quandles are  left-distributive (by idempotency). 
It is proved in \cite{CEHSY}
that  any left-distributive connected quandle is Latin.
This 
 implies, by Toyoda's theorem, that 
every medial connected quandle is Alexander and Latin.
The smallest Latin quandles that are not left-distributive are the Galkin
quandles of order $15$.  It is known that the smallest left-distributive Latin quandle that is not Alexander
is of order $81$. This is due to V. D. Belousov.

\section{Quandle cohomology and cocycle invariant of knots}

In the classical theory of knots and links in 3-space, one utilizes projections of knots and links and applies to them the Reidemeister moves, a sequence of which will take one from any one projection of a given knot or link to any other projection of that knot or link. The Reidemeister moves have played an essential role in the development of a wide variety of invariants for knots and links, since any quantity that remains unchanged by the three moves is an invariant for knots and links. In 1999, Carter et al. \cite{CJKLS} used quandle cohomology  to define combinatorial "state-sum'' invariants for classical knots and knotted surfaces called quandle cocycle invariant (see definition below).  Here we mention some interesting results on surfaces in 4-space they obtained: (1) constructing an example of a sphere that is knotted in 4-dimensional space \cite{CJKLS},  (2) giving obstructions to ribbon concordance for knotted surfaces \cite{CSS}, and (3) detecting non-invertibility of knotted surfaces \cite{CJKLS}.  This was extended to some other examples \cite{CES} and {\cite{CEGS}.\\ In order to define quandle homology and the cocycle knot invariant we need to define coloring of knots by a quandle.   A  {\it coloring}   
of an oriented  classical knot $K$ is a
function ${\mathcal C} : R \rightarrow X$, where $X$ is a fixed 
quandle
and $R$ is the set of over-arcs in a fixed diagram of $K$,
satisfying the  condition
depicted 
in the top
of Figure~\ref{2cocy}. 
This definition of colorings on knot diagrams has been known, see 
\cite{FR} for example. 
In the bottom of Figure~\ref{2cocy}, the relation between Redemeister type III move
and a quandle axiom (self-distributivity) is indicated. 
In particular, the colors of the bottom right segments before and after
the move correspond to the self-distributivity.   By assigning a weight $\phi(x,y)$ at each crossing of a knot diagram (as in the top Figure ~\ref{2cocy}) we obtain a 2-cocycle condition which can be generalized to a homology of cohomology theory which we describe now.

\bigskip
\begin{figure}\label{2cocy}
\begin{center}
\mbox{
\includegraphics[scale=0.4]{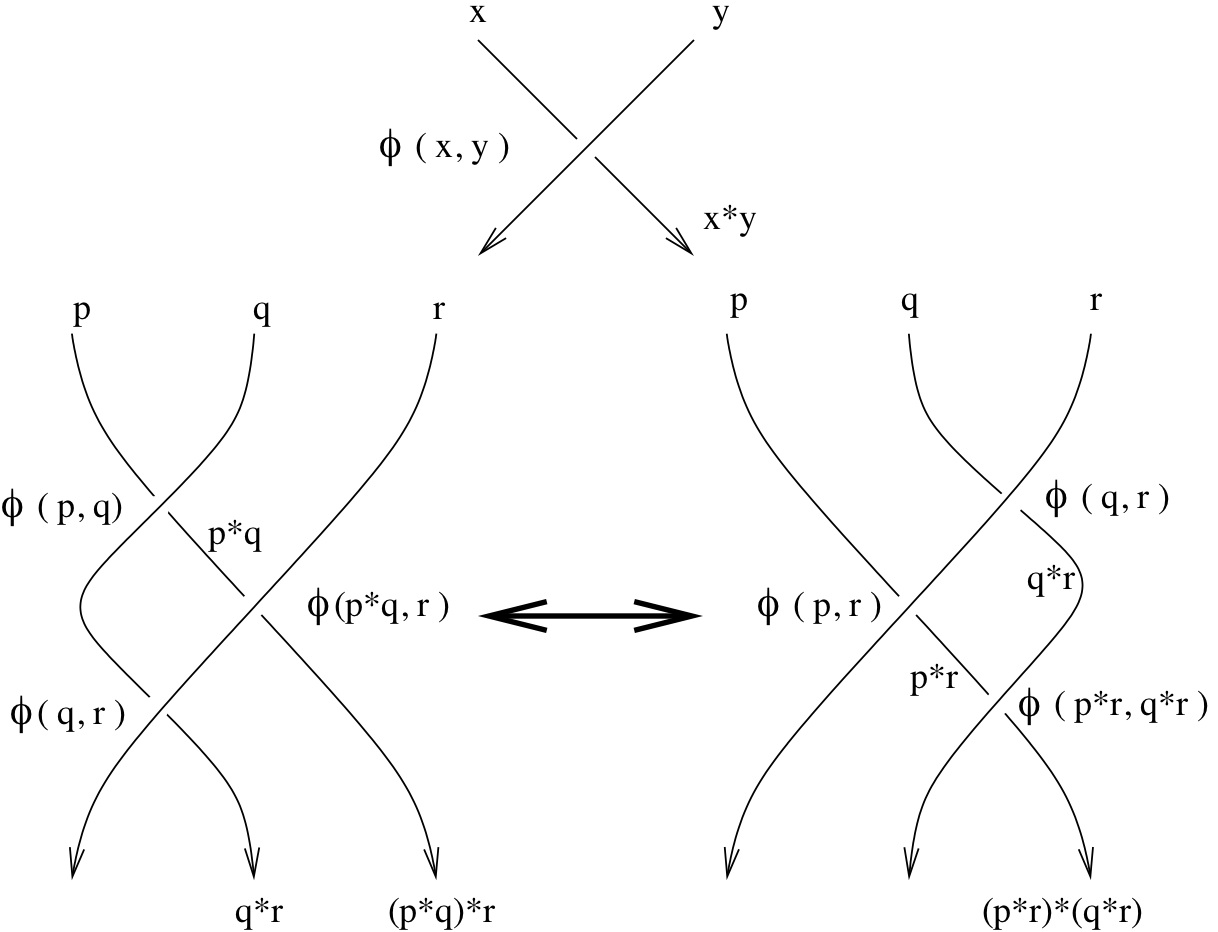}
}
\end{center}
\caption{ 2-cocycle condition coming from Type III move }
\label{qcolor} 
\end{figure}

 Let $C_n(X)$ be the free 
abelian group generated by
$n$-tuples $(x_1, \dots, x_n)$ of elements of a quandle $X$. Define a
homomorphism
$\partial_{n}: C_{n}(X) \to C_{n-1}(X)$ by \begin{eqnarray}
\lefteqn{
\partial_{n}(x_1, x_2, \dots, x_n) } \nonumber \\ && =
\sum_{i=2}^{n} (-1)^{i}\left[ (x_1, x_2, \dots, x_{i-1}, x_{i+1},\dots, x_n) \right.
\nonumber \\
&&
- \left. (x_1 \ast x_i, x_2 \ast x_i, \dots, x_{i-1}\ast x_i, x_{i+1}, \dots, x_n) \right]
\end{eqnarray}
for $n \geq 2$ 
and $\partial_n=0$ for 
$n \leq 1$. 
 Then
$C_\ast(X)
= \{C_n(X), \partial_n \}$ is a chain complex.  The $n$\/th {\it quandle homology group\/}  and the $n$\/th
{\it quandle cohomology group\/ } \cite{CJKLS} of a quandle $X$ with coefficient in a group $A$ can be defined.  One can consider cohomology also and for example: \\ A {\it 2-cocycle} is a  function $\phi: X \times X \rightarrow A$ such that $\phi(x,y)+\phi(x*y,z)=\phi(x,z)+\phi(x*z,y*z),$ and for all $x$, $\phi(x,x)=0$.\\
A {\it 3-cocycle} is a function $\psi: X \times X \times X \rightarrow A$ such that 
$$\psi(x,y,z)+\psi(x,z,w)+\psi(x*z,y*z,w)=\psi(x*y,z,w)+\psi(x*w,y*w,z*w)+\psi(x,y,w),$$
and and for all $x, y $, $\psi(x,x,y)=\psi(x,y,y)=0$.

 Let  ${\cal  C}$ 
denote a coloring of a knot $K$ by a quandle $X$ and choose  a 
quandle 
2-cocycle 
$\phi$,
Then define a {\it (Boltzmann) weight}, $B(\tau, {\cal C})$, at a  crossing $\tau$, by $B(\tau, {\cal C})= \phi(x,y)^{\epsilon (\tau)}$,
where 
$\epsilon (\tau)= 1$ or $-1$, if  the sign of $\tau$ 
is positive or negative, respectively.

The {\it partition function}, or a {\it state-sum}, called quandle cocycle invariant
is the expression 
$$
\Phi_{\phi}(K):=\sum_{{\cal C}}  \prod_{\tau}  B( \tau, {\cal C}).
$$
The product is taken over all crossings of the given diagram,
and the sum is taken over all possible colorings.
The values of the partition function 
are  taken to be in  the group ring ${\Z }[A]$ where $A$ is the coefficient 
group 
written multiplicatively. 
\begin{theorem}  \cite{CJKLS}
The state sum $\Phi_{\phi}(K)$ does not depend on the choice of a diagram of a knot $K$, so that it is a knot invariant.

\end{theorem}

\begin{example} {\rm \cite{CJKS} p 52,
Let $X=\Z_2[T,T^{-1}]/(T^2+T+1)$,  $A=\Z_2$, and cocycle $\Phi=\prod \chi _{(a,b)}$ where $a,\;b \in \{0,1,T+1\}$ and $a \neq b$. \\
For knots $K$ (up to nine crossings, see \cite{CL} for diagrams and other information) the Invariants $\Phi (K)$ are:\\

$\bullet$ $4(1+3T)$ for $3_1, 4_1, 7_2, 7_3, 8_1, 8_4, 8_{11}, 8_{13}, 9_1, 9_6, 9_{12},9_{13},9_{14}, 9_{21}, 9_{23}, 9_{35}, 9_{37},$\\

$\bullet$ $16(1+3T)$ for $8_{18},$ and  $\;9_{40}$\\

$\bullet$ $16$ for $8_5, 8_{10}, 8_{15}, 8_{19}-8_{21},9_{16},9_{22}, 9_{24}, 9_{25}, 9_{28}-9_{30}, 9_{36}, 9_{38}, 9_{39},  9_{41}-9_{45}, 9_{49}$\\

$\bullet$ $4$ otherwise.
}
\end{example}
Generalizations of the cocycle knot invariants have been discovered for example see \cite{CES} and \cite{CEGS}.

\subsection{ Extensions of quandles}
Quandle extension theory was developed in \cite{CENS} by analogy
with group extensions defined for low dimensional group cocycles. Let $X$ be a quandle, $A$ be an abelian group and given a $2$-cocycle $\phi \in Z^2_{\rm Q}(X; A)$, the quandle operation in
extension is defined on
 $E=A \times X$  by
$(a_1, x_1)*(a_2, x_2)=(a_1 + \phi(x_1, x_2),\  x_1 * x_2)$.  
The following lemma is the converse of the fact proved in \cite {CKSdiag} that 
$E(X,A,\phi)$ is a quandle.

\begin{lemma} \label{cocylemma}\cite{CENS}
Let $X$, $E$ be finite quandles, and $A$ be a finite abelian group
written multiplicatively. 
Suppose there exists a bijection 
$f: E \rightarrow A \times X$ with the following property.
There 
exists 
a function $\phi: X \times X \rightarrow A$ such that
for any $e_i \in E$ ($i=1,2$), 
if $f(e_i)=(a_i, x_i)$, then 
$f(e_1 *  e_2) = (a_1 \phi(x_1, x_2) , x_1 * x_2 )$. 
Then $\phi \in Z^2_{\rm Q}(X; A)$. 
\end{lemma}  
The following two theorem produce examples of extensions of quandles.

\begin{theorem} \label{firstextthm} \cite{CENS}
For any  positive integers $q$ and $m$, 
$E=\Z _{q^{m+1}} [T, T^{-1}] /  (T -1 +q) $ is an abelian extension 
$E=E( \Z _{q^{m}} [T, T^{-1}] /  (T -1 +q) ,  \Z_q,  \phi)$
of $X= \Z _{q^{m}} [T, T^{-1}] /  (T -1 +q)$ for some cocycle
 $\phi \in  Z^2_{\rm Q}( X;  \Z_q)$. 
\end{theorem}

\begin{theorem} \label{secondextthm}  \cite{CENS}
For any positive integer $q$ and $m$, the quandle 
$E=\Z_q [T, T^{-1} ] / (1-T)^{m+1} $ is an abelian extension
of $X=\Z_q [T, T^{-1} ] / (1-T)^{m} $ over $\Z_q$: 
$E=E(X, \Z_q, \phi)$, for some $\phi \in Z^2_{\rm Q}(X; \Z_q)$. 
\end{theorem}
Below are some explicit examples of extensions.

\begin{example} \label{secondextthm}
{\rm \cite{CENS}
For any positive integer $q$ and $m$, the quandle 
$E=\Z_q [T, T^{-1} ] / (1-T)^{m+1} $ is an abelian extension
of $X=\Z_q [T, T^{-1} ] / (1-T)^{m} $ over $\Z_q$: 
$E=E(X, \Z_q, \phi)$, for some $\phi \in Z^2_{\rm Q}(X; \Z_q)$. 
}
\end{example}

\begin{example} {\rm  \cite{CENS}
Consider the case $q=2$, $m=2$ in Example~\ref{secondextthm}.
In this case 
$$\Z_4[T, T^{-1}]/ (T+1)=R_4, \quad \mbox{and} $$ 
 $$\Z_8[T, T^{-1}]/ (T+1)=R_8=E(R_4, \Z_2, \phi)$$ for some 
$\phi \in Z^2_{\rm Q}(R_4; \Z_2)$.
We  obtain an explicit formula for this cocycle
$\phi$ by computation: 
$$\phi = \chi_{0,2} + \chi_{0,3} + \chi_{1,0} + \chi_{1,3}
+ \chi_{2,0} + \chi_{2,3} + \chi_{3,0} + \chi_{3,1}   , $$
where 
$$\chi_{a,b} (x,y) = \left\{ \begin{array}{ll} 1 & {\mbox{\rm if }} \
(x,y)=(a,b), \\
0 & {\mbox{\rm if }} \
(x,y)\not=(a,b) \end{array}\right.$$
denotes the characteristic function.

} \end{example}
Other extensions of quandles have been considered by some authors see for example in \cite{AG} where a more general homology theory is developed and in \cite{ Eiserman} where algebraic covering theory of quandle is established.\\
{\bf Dynamical cocycles}  \cite{AG} 
Let $X$ be a quandle and $S$ be a non-empty set. 
Let $\alpha: X \times X \rightarrow 
\mbox{\rm Fun}(S \times S, S)=S^{S \times S}$ be a function,
so that for $x,y  \in X$ and $a, b \in S$ we have 
$\alpha_{x,y} (a,b) \in S$. 

Then it is checked by computations 
that $S \times X$ is a quandle by the operation
$(a, x)*(b, y)=(\alpha_{x, y}(a,b )  , x*y )$,
where $ x*y$ denotes the quandle 
product in $X$, 
if and only if $\alpha$ satisfies the following conditions:

\begin{enumerate}
\setlength{\itemsep}{-3pt}

\item $\quad\alpha_{x, x} (a,a)= a$ for all $x \in X$ and 
$a \in S$; 
\item $\quad\alpha_{x,y} (-,b): S \rightarrow S$ is a bijection for 
all $x, y \in X$ and for all $b \in S$;
\item $\quad 
\alpha_{x * y, z} ( \alpha_{x , y}(a, b), c)
=\alpha_{x * z, y*z}(\alpha_{x,z}(a,c), \alpha_{y,z}(b,c) )$
for all $x, y, z \in X$ and $a,b,c \in S$. 

\end{enumerate}

Such a function $\alpha$ is called a {\it dynamical quandle cocycle}
\cite{AG}.
The quandle constructed above is denoted by $S \times_{\alpha} X$, 
and is called the {\it extension} of $X$ by a dynamical cocycle $\alpha$.
The construction is general, as Andruskiewitsch and Gra\~{n}a
 show:

\begin{lemma}{ \cite{AG}} \label{AGlemma} 
Let $p: Y \rightarrow X$ be a surjective quandle homomorphism 
between finite quandles 
such that the cardinality of $p^{-1}(x)$ is a constant for all $x \in X$. 
Then $Y$ is isomorphic to an extension  $S \times_{\alpha} X$ of $X$ 
by some dynamical cocycle on 
the
set $S$ 
such that $|S|=|p^{-1}(x)|$.
\end{lemma}

\noindent
{\large\bf Acknowledgments} \  I would like to thank Scott Carter, Edwin Clark and Masahico Saito for fruitful suggestions.

\end{document}